\documentclass[12pt]{article}
\usepackage{amsmath,  amssymb,amsthm}
\newtheorem{theorem}{Theorem}

\newtheorem{lemma}[theorem]{Lemma}

\theoremstyle{definition}
\newtheorem{remark}[theorem]{Remark}

\begin{document}

\title{Chern classes in Alexander-Spanier cohomology}
\author{Alexander Gorokhovsky\\Department of Mathematics,\\
  Ohio State
University,\\
 Columbus, OH 43210\\
sasha@math.ohio-state.edu}
\maketitle

\begin{abstract}
In this article we construct explicit cocycles in the Alexander-Spanier
cohomological complex, representing
 the Chern character of an element in K-theory
\end{abstract}

\section{Introduction. }
 In  recent years   explicit formulas for
 cocycles  representing
 the Chern character of an element in $K$-homology
were obtained in several papers ([CST], [MW]).  This led  in particular
to a solution of the  problem of finding explicit formulas for the
rational Pontrjagin classes of  topological manifolds. 
One of the  distinctive features of these approaches,  based on [CM],  is the
use of the Alexander-Spanier (co)homology,  so the cocycles obtained live in
the Alexander-Spanier homological complex. 

In connection with this  a ``dual'' problem naturally arises, that
of finding  explicit
formulas, representing the Chern character of an element in  $K$-theory in 
the Alexander-Spanier cohomology. In this paper we construct such  cocycles. 
In conjunction with the formulas from [CST], [MW] this allows to write down
explicitly index pairing in terms  of the Alexander-Spanier (co)homology.

Let us briefly outline the main steps of our construction. 
Alexander-Spanier $n$-cochains are continuous functions in the
neighborhood of the diagonal in $X^{n+1}$. 
Now,  for example,  if we consider an element in $K^0(X)$  represented
by an idempotent $e(x)\in M_N(C(X))$, our cocycle computed
at the point $(x_0, x_1, \dots, x_n)\in X^{n+1}$
is an integral of the differential form
representing the Chern character of the universal
bundle over the Grassmanian over the canonical
simplex with the vertices at $e(x_0)$, $e(x_1)$, \dots, $e(x_n)$. 
It follows easily from  Stokes' theorem
that the cochain constructed in this way is a cocycle. 
To verify that our cocycle represents the Chern character components in the smooth case
we compute its image under the canonical projection  to
de Rham cohomology. This projection effectively 
amounts to averaging our form over  smaller and smaller simplices, so
in the limit one should get the original form. 
Finally,  for the case of a  general compact topological space
the result follows by the argument
which uses functoriality and homotopy invariance of classes
of our cocycles.

The paper is organized as follows : in  section \ref{prel}
we recall the main facts about the Alexander-Spanier complex, 
and in  sections \ref{odd} and \ref{even} we construct
cocycles representing the Chern character of elements
in $K^{-1}(X)$ and $K^0(X)$ respectively.

\section{Preliminaries. }\label{prel}
We will describe a version of the Alexander-Spanier cohomology (with
complex coefficients),  following [CM],  [MW].  All the proofs can be found
in these papers.  Let us recall the
main definitions and facts concerning this cohomology. 
Let $X$ be a compact separable topological space. 
Let $Cov^f(X)$ denote the set of all finite open coverings of $X$,  and let $\mathfrak{U}\in 
Cov^f(X)$.  Let $\mathfrak{U} ^{n}$ denote the neighborhood
of the diagonal in $X^n$ given by $\cup_{U\in \mathfrak{U}} U^n$.  Then Alexander-Spanier
$n$-cocycles (corresponding to $\mathfrak{U}$) are continuous functions on $\mathfrak{U} ^{n+1}$. 
 The space of $n$-cocycles is denoted by $C^n(X, \mathfrak{U})$.  The boundary operator
$\partial:C^n(X, \mathfrak{U})\to C^{n+1}(X, \mathfrak{U})$ is defined by the formula
(here $\phi \in C^n(X, \mathfrak{U})$)  
\begin{equation}\label{defb}
\partial \phi(x_0, x_1, \dots, x_{n+1})=\sum_{j=0}^{n+1}(-1)^j \phi (x_0, \dots, \hat x_j, \dots, x_{n+1})
\end{equation}
Then $\partial^2=0$ and the cohomology of the complex $(C^n(X, \mathfrak{U}), \partial)$ is called
Alexander-Spanier cohomology of the covering $\mathfrak{U}$; it is denoted $H^n(X, \mathfrak{U})$.
If $\mathfrak{V}$ is a refinement of $\mathfrak{U}$,  there is an obvious map (restriction)
from $C^n(X, \mathfrak{U})$ to $C^n(X, \mathfrak{U})$,  which commutes with the boundary
from (\ref{defb}) ,  and hence defines a map $ H^n(X, \mathfrak{U}) \to H^n(X, \mathfrak{V}) $. 
Then Alexander-Spanier cohomology is defined as a direct limit over finite covers:
\[
H^n(X)=\varinjlim H^n(X, \mathfrak{U})
\]
 The cup-product on the Alexander-Spanier complex is given by the formula 
 \[
 \phi \cup \psi (x_0, \dots, x_{m+n})=\phi (x_0, \dots, x_n)\psi (x_n, \dots, x_{m+n}) 
\]
where $\phi$ and $\psi$ are $n$ and $m$ cochains respectively. 

We will also describe, following [MW],  the dual theory --- Alexander-Spanier homology. 
Again,  we fix a cover $\mathfrak{U}$ of $X$.  Let $C_n(X, \mathfrak{U})$ be the space
of measures with (compact) support in $\mathfrak{U}^{n+1}$ --- chains for the Alexander-Spanier
homology complex.  This is the dual space for $C^n(X, \mathfrak{U})$ endowed with a topology
of uniform convergence on the compacts. 
The boundary operator,  which we also denote $\partial$, is defined as transposed to the operator
given by (\ref{defb}):
\[
\partial \mu(\phi):=\mu(\partial \phi)
\]
Here $\mu \in C_n(X, \mathfrak{U})$,  $ \phi \in C^n(X, \mathfrak{U})$.
The homology of $(C_n(X, \mathfrak{U}), \partial)$ is called
Alexander-Spanier homology of $\mathfrak{U}$ and 
is denoted $H_n(X, \mathfrak{U})$. 
Now,  if $\mathfrak{V}$ is a refinement of $\mathfrak{U}$,  there is a map of complexes
$C_*(X, \mathfrak{V}) \to C_*(X, \mathfrak{U})$ (extension by zero). It induces a map 
$ H_n(X, \mathfrak{V}) \to H_n(X, \mathfrak{U}) $. Alexander-Spanier homology of $X$ is then defined
as an inverse limit over finite covers: 

\[
H_n(X)=\varprojlim H_n(X, \mathfrak{U})
\]

It can be shown that Alexander-Spanier (co)homology is isomorphic to the \v{C}ech (co)homology, 
and hence to the de Rham (co)homology if $X$ is a smooth manifold. In the latter case one can
consider the smooth Alexander-Spanier complex,  where $n$-cochains are given by  smooth functions
in the neighborhood of the diagonal in $X^{n+1}$,  and the coboundary  operator is given by (\ref{defb}).
We will now describe a canonical morphism $\lambda$ from the smooth Alexander-Spanier complex
to the de Rham cohomological complex,  which induces an isomorphism in cohomology. 
Let $\phi$ be a smooth Alexander-Spanier $n$-cocycle. Let $x_0 \in X$,  and $v_1, \dots, v_n \in
T_{x_0} X$.  Consider any $n$ curves $x_j(\epsilon)$,  $j=1, \dots, n$ with the properties
\begin{align}
x_j(0)=x_0\notag\\
\frac{d}{d\epsilon} x_j |_{\epsilon=0}=v_j\notag
\end{align}
Then the differential $n$-form $\lambda (\phi)$ is defined by the equation:
\begin{multline}
\lambda (\phi)_{x_0}(v_1, \dots, v_n)\\
=\frac{1}{n!} \sum_{\sigma \in S_n}
sgn(\sigma) \frac{\partial}{\partial \epsilon_1} \dots \frac{\partial}{\partial \epsilon_n}
\phi \left(x_0,  x_{\sigma(1)}(\epsilon_1),  \dots,  x_{\sigma(n)}(\epsilon_n)\right) |_{\epsilon_j=0}
\end{multline}

\section{The odd case. }\label{odd}
In this section we will give formulas for the Alexander-Spanier
cochain
describing the odd Chern character
$Ch:K^{-1}(X)\rightarrow H^{odd}(X)$,  where $X$ is a compact topological space.
Corresponding formulas for the de Rham cohomology in the smooth case are
well known; a good exposition can be found in the first section of [Ge].
An element in $K^{-1}$ can be represented by a continuous map
$U$ from $X$ to $\mathbf{U} (N )$ --
the group of unitary $N\times N$ matrices,  for some $N$ (note that since we are working with 
complex
coefficients,  $K^{-1}(X)=K^1(X)$). 

Now,  let $x_0, x_1, \dots, x_n$ be $n+1$ sufficiently close points in $X$,  so that
$\|U(x_k)-U(x_l)\|<\rho<1$ for some $\rho$, and where $\| \cdot \|$
is the operator norm. 

Let $t_1, t_2, \dots, t_n$ be nonnegative numbers with $t_1+t_2 \dots +t_n\leq 1$. 
Put $t_0=1-\sum_{j=0}^n t_j$. 
Define
\begin{equation}\label{defu}
U(t_1, t_2, \dots, t_n)=\sum_{j=0}^n t_jU(x_j)=U(x_0)+\delta
\end{equation}
where
\begin{align}\label{defdu}
\delta=\sum_{j=0}^n t_j(U(x_j)-U(x_0))=\sum_{j=0}^nt_j \delta_j\\
\delta_j=U(x_j)-U(x_0)
\end{align}
Since $\|\delta \|<1 $, $U(t_1, t_2, \dots, t_n)$ is always invertible.  

We can now consider on the $n$-simplex
\[ \Delta^n=
\{t_j, j=1, \dots, n|\sum_{j=1}^n t_j\leq 1\}=
\{t_j, j=0, \dots, n|\sum_{j=0}^n t_j=1\} \]
 a
 matrix-valued 1-form
$U(t_1, \dots, t_n)^{-1}dU(t_1, \dots, t_n)$, 
where
\[
dU(t_1, \dots, t_n)=\sum_{j=1}^n \frac{\partial U(t_1, \dots, t_n)}{\partial t_j}dt_j
\]
This form can be rewritten as
\begin{multline}\label{frmu}
U(t_1, t_2, \dots, t_n)^{-1}dU(t_1, t_2, \dots, t_n)=(U(x_0)+\delta)^{-1}d(U(x_0)+\delta)\\
=\sum_{j=1}^n \sum_{k=0}^{\infty}(-1)^k(U(x_0)^{-1}\delta)^k U(x_0)^{-1}\delta_jdt_j
\end{multline}
From this one can construct for $n$ odd an Alexander-Spanier cochain
\begin{equation}\label{defcu}
Ch^n(x_0, \dots, x_n)=c_n\int_{\Delta^n}
Tr\
( U(t_1, t_2, \dots, t_n)^{-1}dU(t_1, t_2, \dots, t_n))^n
\end{equation}

Here  the $n$-th power is taken with respect to the exterior product
and $c_n=\frac{(-1)^{(n-1)/2}((n-1)/2)!}{(2\pi i)^{(n+1)/2}n!}$.

For example, for $n=1$, $\delta=t_1\delta_1=t_1(U(x_1)-U(x_0))$, 
\begin{multline}
U(t_1)^{-1}dU(t_1)=\sum_{k=0}^{\infty}(-1)^k(U(x_0)^{-1}\delta)^k U(x_0)^{-1}\delta_1dt_1\\
=\sum_{k=0}^{\infty}(-1)^k(U(x_0)^{-1}\delta_1)^{k+1}t_1^kdt_1
\end{multline}
by (\ref{frmu}),  and
\begin{multline}
Ch^1(x_0, x_1)
=\frac{1}{2\pi i}\int_0^1Tr\sum_{k=0}^{\infty}(-1)^k(U(x_0)^{-1}\delta_1)^{k+1}t_1^kdt_1\\
=\frac{1}{2\pi i}Tr\sum_{k=0}^{\infty}\frac{(-1)^k(U(x_0)^{-1}\delta_1)^{k+1}}{k+1}=
\frac{1}{2\pi i}Tr\log(1+U(x_0)^{-1}\delta_1)\\
=\frac{1}{2\pi i}\log \det \left((1+U(x_0)^{-1}\delta_1\right)=
\frac{1}{2\pi i}\log \det\left( U(x_0)^{-1}U(x_1)\right)
\end{multline}
\begin{remark}
Note that this function is well-defined and continuous in the
neighborhood of the diagonal
given by condition $\|\delta_1\|<\rho<1$,  while  the function
$\log \det U(x_0)$
is not well defined. This explains why we \emph{cannot} write that
\begin{multline*}
Ch^1(x_0,x_1)=\frac{1}{2\pi i}
\left(-\log \det U(x_0)+\log \det U(x_1) \right)\\=
\partial \left(\frac{1}{2\pi i} (\log \det U)(x_0,x_1) \right)
\end{multline*}
\end{remark}
We can now formulate the main result of this section:
\begin{theorem}\label{thmu}
Cochains $Ch^n$,  defined by (\ref{defcu}),  give components for the Chern character of $U(x)$. 
\end{theorem}
To prove this theorem we need several lemmas. 
\begin{lemma}\label{closu}
Formula (\ref{defcu}) defines an Alexander-Spanier cocycle. 
\end{lemma}
\begin{proof}
First,  notice that $Ch^n$ is clearly a continuous function ( in the neighborhood of the diagonal).  
Now suppose that we are given $n+2$ points $x_0,  \dots,  x_{n+1}$.
Similarly to (\ref{defu}) we can construct an operator over
$\Delta^{n+1}=\{t_j,  j=0,  \dots,  n+1 \mid \sum_{j=0}^{n+1} t_j =1,  t_j \geq 0\}$

\begin{equation}
V(t_1, t_2, \dots, t_{n+1})=\sum_{j=0}^{n+1} t_jU(x_j)
\end{equation}

We can consider for any $n$ a form $Tr(V^{-1}dV)^n$ on $\Delta^{n+1}$. 
This form is 0 for even $n$ --- indeed, 
\begin{multline}
Tr(V^{-1}dV)^n=TrV^{-1}dV(V^{-1}dV)^{n-1}\\
=(-1)^{n-1}Tr(V^{-1}dV)^{n-1}V^{-1}dV=-Tr(V^{-1}dV)^n \notag
\end{multline}
For $n$ odd this form is closed,  since, using the 
equality $d(V^{-1})=-V^{-1}dVV^{-1}$ we get:
\[
dTr(V^{-1}dV)^n=-Tr(V^{-1}dV)^{n+1}=0
\]
by the previous computation. 

Now, the restriction of $V$ to the face given by equation $t_j=0$ coincides (after an obvious
renumeration of variables) with the operators
constructed by the formula (\ref{defu}) from the points $x_0,  \dots,  \hat x_j,  \dots,  x_{n+1}$.
From this we conclude that
\begin{multline}
\partial Ch^n(x_0,  \dots,  x_{n+1})=\sum_{j=0}^{n+1} (-1)^jCh^n(x_0,  \dots,  \hat x_j,  \dots,  x_{n+1})\\
=\frac{(-1)^{(n-1)/2}((n-1)/2)!}{(2\pi i)^{(n+1)/2}n!}\int_{\partial \Delta^{n+1}} Tr(V^{-1}dV)^n\\=
\frac{(-1)^{(n-1)/2}((n-1)/2)!}{(2\pi i)^{(n+1)/2}n!}\int_{\Delta^{n+1}}dTr(V^{-1}dV)^n=0
\end{multline}
by  Stokes' theorem. 
\end{proof}

\begin{lemma}\label{indu}
The cohomology class of the cocycle defined by (\ref{defcu}) depends only on 
the  class $[U] \in K^{-1}(X)$ and
not on the particular representative $U$.
\end{lemma}
\begin{proof}
Let $U_0(x)$,  $U_1(x)$ be two maps
into $\mathbf{U} (N )$ representing the same element in $K^{-1}(X)$. 
By taking $N$ big enough one can suppose that $U_0(x)$ and  $U_1(x)$ are homotopic
via piecewise-smooth family of maps $U_{\tau}(x)$, $0\leq \tau \leq 1$.
 Using the formula (\ref{defu}) we define
homotopy $U_{\tau}(t_1, \dots, t_n)$ between $U_0(t_1, \dots, t_n)$ and $U_1(t_1, \dots, t_n)$. 

We will now show that cocycles corresponding to
the smoothly homotopic unitaries differ by a coboundary. 
Define $Ch^n_{\tau}$ by the equation (\ref{defcu}). 
Then (we write just $U_{\tau}$ for $U_{\tau}(t_1, \dots, t_n)$):  
\begin{equation}\label{dcdt}
\frac{d}{d\tau}Ch^n_{\tau}=\frac{(-1)^{(n-1)/2}((n-1)/2)!}{(2\pi i)^{(n+1)/2}n!}\int_{\Delta^n}\frac{d}{d\tau}
Tr(U_{\tau}^{-1}dU_{\tau})^n
\end{equation}
But
\[
\frac{d}{d\tau}Tr(U_{\tau}^{-1}dU_{\tau})^n
=n\left(d TrU_{\tau}^{-1}\frac{dU}{d\tau}(U_{\tau}^{-1}dU_{\tau})^{n-1}\right)
\]
Indeed, 
\begin{multline}
\frac{d}{d\tau}Tr(U_{\tau}^{-1}dU_{\tau})^n=
nTr\left(\frac{d}{d\tau}(U_{\tau}^{-1}dU_{\tau})\right)(U_{\tau}^{-1}dU_{\tau})^{n-1}\\
=nTr\left(U_{\tau}^{-1}\frac{dU}{d\tau}U_{\tau}^{-1}dU_{\tau}(U_{\tau}^{-1}dU_{\tau})^{n-1}\right)+
nTr\left(U_{\tau}^{-1}d(\frac{dU}{d\tau})(U_{\tau}^{-1}dU_{\tau})^{n-1} \right)\\
=nTr\left(U_{\tau}^{-1}dU_{\tau}U_{\tau}^{-1}\frac{dU}{d\tau}(U_{\tau}^{-1}dU_{\tau})^{n-1}\right)+
nTr\left(U_{\tau}^{-1}d(\frac{dU}{d\tau})(U_{\tau}^{-1}dU_{\tau})^{n-1} \right)
\end{multline}
and
\begin{multline}
d\ TrU_{\tau}^{-1}\frac{dU}{d\tau}(U_{\tau}^{-1}dU_{\tau})^{n-1}=
Tr\ d(U_{\tau}^{-1})\frac{dU}{d\tau}(U_{\tau}^{-1}dU_{\tau})^{n-1}\\
+TrU_{\tau}^{-1}d(\frac{dU}{d\tau})(U_{\tau}^{-1}dU_{\tau})^{n-1}
+TrU_{\tau}^{-1}\frac{dU}{d\tau}d(U_{\tau}^{-1}dU_{\tau})^{n-1}\\
=Tr\ U_{\tau}^{-1}( dU_{\tau})U_{\tau}^{-1}\frac{dU}{d\tau}(U_{\tau}^{-1}dU_{\tau})^{n-1}+
TrU_{\tau}^{-1}d(\frac{dU}{d\tau})(U_{\tau}^{-1}dU_{\tau})^{n-1}
\end{multline}
since $d(U_{\tau}^{-1}dU_{\tau})^{n-1}=0$ for odd $n$. 

Hence, from (\ref{dcdt})

\begin{equation}\label{derc}
\frac{d}{d\tau}Ch^n_{\tau}=\int_{\partial \Delta^n}\frac{(-1)^{(n-1)/2}((n-1)/2)!}{(2\pi i)^{(n+1)/2}n!}
n\left( TrU_{\tau}^{-1}\frac{dU}{d\tau}(U_{\tau}^{-1}dU_{\tau})^{n-1}\right) 
\end{equation}
 Now notice that since the restriction of $U_{\tau}$ to the $j$-th face of the $\Delta^n$ given
 by equation $t_j=0$ depends only on $U(x_k), k\neq j$, and not on $U(x_j)$, 
 we can define an Alexander-Spanier $n-1$ cochain
 \[
 T_{\tau}(x_0, \dots, x_{n-1})=\frac{(-1)^{(n-1)/2}((n-1)/2)!}{(2\pi i)^{(n+1)/2}n!}
n\int_{\Delta^{n-1}}
TrU_{\tau}^{-1}\frac{dU}{d\tau}(U_{\tau}^{-1}dU_{\tau})^{n-1}
\]
where $\Delta^{n-1}$=$\{(t_1, \dots,  t_{n-1})\mid t_j\geq 0, \sum_{j=1}^{n-1} t_j \leq 1\}=$
$\{(t_1, \dots,  t_{n})\in \Delta^n\mid t_n=0\}$. 
Then according to (\ref{derc})
\[
 \frac{d}{d\tau}Ch^n_{\tau}=\partial T_{\tau}
\]
and
\[
Ch^n_1-Ch^n_0=\partial \int_0^1 T_{\tau}d{\tau}
\]
which proves the Lemma.
\end{proof}
 
\begin{lemma}\label{smthu}
Let $X$ be a compact smooth manifold and $U:X\rightarrow \mathbf{U}(N )$ be a smooth map.
Then $Ch^n$ given by the formula (\ref{defcu}) represents the $n$-th component of the Chern
character of $ [U]\in K^{-1}(X)$. 
\end{lemma}
\begin{proof} It is known that in the situation described in the Lemma
the differential form
\begin{equation}\label{defdiffcu}
\Omega=\frac{(-1)^{(n-1)/2}((n-1)/2)!}{(2\pi i)^{(n+1)/2}n!}Tr(U^{-1}dU)^n
\end{equation}
 represents
the $n$-th component of the Chern character ( here  the differential is taken
with respect to the $x $ variable). We will show that the canonical map
from Alexander-Spanier to de Rham complex maps cocycle given by (\ref{defcu})
in the differential form given by (\ref{defdiffcu}). 
  We consider a point $x_0$ and $n$ curves $x_j(\epsilon_j),  j=1,  \dots,  n$ such that
$x_j(0)=x_0$ ; let the tangent vector
to $x_j$ at $x_0$ be  $v_j$. 
From (\ref{defdu}) we  have easily
\begin{align*}
\delta|_{\epsilon_k=0}=0\\
\delta_j|_{\epsilon_k=0}=0
\end{align*}
Also
\[
\frac {\partial \delta} {\partial \epsilon_k} |_{\epsilon_j=0}=t_k(\mathcal{L}_{v_k}U)(x_0)
\]
and \[
\frac {\partial \delta_l}{\partial \epsilon_k} |_{\epsilon_j=0}=
\begin{cases}
(\mathcal{L}_{v_k}U)(x_0), &\text{if $k=l$}\\
  0, &\text{otherwise}
\end{cases}
\] 
Here $\mathcal{L}_{v_k}$ denotes directional derivative. 
From this by  differentiating (\ref{frmu}) one gets
\[
\frac{\partial U(t_1, t_2, \dots, t_n)^{-1}dU(t_1, t_2, \dots, t_n) }
{\partial \epsilon_k}|_{\epsilon_j=0}=U^{-1}(x_0)(\mathcal{L}_{v_k}U)(x_0)dt_k
\]

We now compute with $c_n=\frac{(-1)^{(n-1)/2}((n-1)/2)!}{(2\pi i)^{(n+1)/2}n!}$
\begin{multline}
\frac{\partial}{\partial \epsilon_1} \dots \frac{\partial}{\partial \epsilon_n}
Ch^n\left( x_0,  x_1(\epsilon_1),  \dots,  x_n(\epsilon_n) \right) |_{\epsilon_j=0}\\
=c_n\int_{\Delta^n}Tr
\frac{\partial}{\partial \epsilon_1} \dots \frac{\partial}{\partial \epsilon_n}
(U(t_1, t_2, \dots, t_n)^{-1}dU(t_1, t_2, \dots, t_n))^n |_{\epsilon_j=0}\\
=\sum_{\sigma\in S_n} sgn(\sigma)
c_n\int_{\Delta^n}Tr
U^{-1}(x_0)(\mathcal{L}_{v_{\sigma(1)}}U)(x_0) \dots
U^{-1}(x_0)(\mathcal{L}_{v_{\sigma(n)}}U)(x_0)dt_1\dots dt_n\\
=c_n\sum_{\sigma\in S_n} sgn(\sigma)
\frac{1}{n!}Tr
U^{-1}(x_0)(\mathcal{L}_{v_{\sigma(1)}}U)(x_0) \dots
U^{-1}(x_0)(\mathcal{L}_{v_{\sigma(n)}}U)(x_0)\\
=\Omega_{x_0}(v_1, \dots, v_n)
\end{multline}
The result is already antisymmetric in $v_1, \dots, v_n$,
and the assertion follows. 
\end{proof}

Now we can prove our Theorem \ref{thmu}.

\begin{proof}[Proof of  Theorem \ref{thmu}. ]
 When $X$ and $U$ are smooth the result was just proved. 
 For the general $X$ and $U$, since any element in $K^{-1}(X)$ can
 be pulled back from some smooth manifold,  and our construction is manifestly
 functorial with respect to pull-backs of unitaries,  we find that for any class
 in $K^{-1}(X)$ for \emph{some} representative our formula represents the components
 of the Chern character. But then  by the Lemma \ref{indu} this is true for \emph{all}
representatives.  
\end{proof}

\section{The even case. }\label{even}

Let $X$ be a compact topological space and $E\rightarrow X$ be a complex
vector bundle.  Suppose that $E$ is embedded as a subbundle into a trivial
bundle with total space $X\times \mathbb{C}^N$.  This allows us to represent
$E$ by a self-adjoint projector $e(x)$ in $M_N(C(X))$ --- the algebra of $N\times N$
matrices of continuous functions on $X$.   
We will construct here an Alexander-Spanier cochain
representing the Chern character of $E$.  

Let $x_0,  x_1,  \dots,  x_n$ be $n+1$ points in $X$  sufficiently close to each other:
namely, we require that $\|e(x_i)-e(x_j)\|<\rho<1/2$ (with respect to the usual operator
 norm).  We will now define the function $Ch^n(E)$,  representing the $n$-th component of the
Chern character of $E$ ($n$ is even here).  Let $t_1,  t_2,  \dots,  t_n$ be positive numbers with
$t_1+t_2+\dots+t_n \leq 1$.  Put $t_0=1-\sum_{i=1}^n t_i$.  Consider the 
operator
\[
a(t_1,  \dots,  t_n)=\sum_{i=0}^n t_ie(x_i)=e(x_0)+\delta
\]
 where
 \begin{equation}\label{defd}
 \delta=\sum_{i=1}^n t_i(e(x_i)-e(x_0)) 
 \end{equation}
Since $\|a(t_1,  \dots,  t_n)-e(x_0)\|\leq \rho$, the  spectrum of $a(t_0,  \dots,  t_n)$ is contained
in the union of 2 discs of radius $\rho<1/2$ and centers at 0 and 1.  
Let $e(t_1,  \dots,  t_n)$ be the  spectral projector on  the part of the spectrum inside the disc
$|\lambda-1|<1/2$:
\begin{equation}\label{defe}
 e(t_1,  \dots,  t_n)=\frac{1}{2\pi i}\int_{|\lambda-1|=1/2} (\lambda-a(t_1,  \dots,  t_n))^{-1}d\lambda
\end{equation}
Notice that $ e(t_1,  \dots,  t_n)$ depends smoothly on $t_1,  \dots,  t_n$,  since
$a(t_1,  \dots,  t_n)$ depends smoothly on them and its spectrum does not intersect the contour
of integration.  

Now  we can rewrite $ e(t_1,  \dots,  t_n)$ in yet another form.

\begin{multline}
 e(t_1,  \dots,  t_n)=\frac{1}{2\pi i}\int_{|\lambda-1|=1/2} (\lambda-a(t_1,  \dots,  t_n))^{-1}d\lambda\\
=\frac{1}{2\pi i}\int_{|\lambda-1|=1/2} (\lambda-(e(x_0)+\delta))^{-1}d\lambda\\
=\frac{1}{2\pi i}\int_{|\lambda-1|=1/2} (1-(\lambda-e(x_0))^{-1}\delta)^{-1}(\lambda-e(x_0))^{-1}d\lambda\\
=\frac{1}{2\pi i}\sum_{k=0}^\infty 
\int_{|\lambda-1|=1/2} ((\lambda-e(x_0))^{-1}\delta)^k(\lambda-e(x_0))^{-1}d\lambda
\end{multline}
We now use the identity
\[
(\lambda-e(x_0))^{-1}=\frac{e(x_0)}{\lambda-1}+\frac{1-e(x_0)}{\lambda}
\]
to rewrite the $k$-th term:
\begin{multline}
\frac{1}{2\pi i}\int_{|\lambda-1|=1/2} ((\lambda-e(x_0))^{-1}\delta)^k(\lambda-e(x_0))^{-1}d\lambda\\
= \frac{1}{2\pi i}\int_{|\lambda-1|=1/2}
\left(\left(\frac{e(x_0)}{\lambda-1}+\frac{1-e(x_0)}{\lambda}\right)\delta\right)^k
\left(\frac{e(x_0)}{\lambda-1}+\frac{1-e(x_0)}{\lambda}\right)d\lambda\\
=\frac{1}{2\pi i}\int_{|\lambda-1|=1/2} \sum_{m+l=k+1}
\frac{1}{(\lambda-1)^{m}\lambda ^{l}}( b_0\delta b_1\delta \dots \delta b_k) d\lambda
\end{multline}
Here each $b_j$ equals either $e(x_0)$ or $(1-e(x_0))$; each such monomial
contains m factors of the first type and l factors of the second.
We see that the expression under the integral has the only pole inside the contour of integration
at the point 1,  and the residue at this point can be explicitly computed.  
Indeed,  by the binomial formula,  $(\lambda-1)^{-m}\lambda^{-l}=\sum_{p=0}^\infty
\binom{-l}{p} (\lambda-1)^{-m} (\lambda-1)^{-p}$,  and hence $res |_{\lambda=1}
(\lambda-1)^{-m}\lambda^{-l}=\binom{-l}{m-1}=(-1)^{m-1}\binom{l+m-2}{m-1}=(-1)^{m-1}
\binom{k-1}{m-1}$ ( here we suppose, as usual,
that $\binom{k-1}{-1} =0$ ) ,  and

\begin{multline}
\frac{1}{2\pi i}\int_{|\lambda-1|=1/2} ((\lambda-e(x_0))^{-1}\delta)^k(\lambda-e(x_0))^{-1}d\lambda\\
=\sum (-1)^{m-1}\binom{k-1}{m-1} b_0\delta b_1\delta \dots \delta b_k
\end{multline}
Here in the term with the coefficient $(-1)^{m-1}\binom{k-1}{m -1}$
there are  $m$ factors  $b_j$  equal to $e(x_0)$.  

Hence
\begin{equation}\label{frme}
e(t_1,  t_2,  \dots,  t_n)=\sum_k \sum (-1)^{m-1}
\binom{k-1}{m-1} b_0\delta b_1\delta \dots \delta b_k
\end{equation}

 For $n$ even we define a function
\begin{equation}\label{defc}
Ch^n(x_0,  x_1,  \dots,  x_n)=b_n\int\limits_{\sum_{j=1}^n t_j\leq 1,  t_j \geq 0}
Tr\ e(t_1,  t_2,  \dots,  t_n)(de(t_1,  t_2,  \dots,  t_n))^n
\end{equation}
with $b_n=\frac{(-1)^{n/2}}{(2\pi i)^{n/2}(n/2)!}$.
Here $de(t_1,  t_2,  \dots,  t_n)=\sum_{j=1}^n \frac{\partial e(t_1,  t_2,  \dots,  t_n)}{\partial t_j}dt_j$
--- a matrix-valued
1-form on the $n$-simplex $\Delta^n=\{t_j,  j=1,  \dots,  n \mid \sum_{j=1}^n t_j\leq 1,  t_j \geq 0\}=
\{t_j,  j=0,  \dots,  n \mid \sum_{j=0}^n t_j =1,  t_j \geq 0\}$

We are now going to prove that $Ch^n$ represents the $n$-th component of the Chern character of the
bundle $E$ . 
\begin{theorem}\label{prf}
Let $E$ be a complex vector bundle over a compact topological space $X$,  represented by the idempotent $e(x)$,  
and let $Ch^n$ be the Alexander-Spanier cocycle
defined by (\ref{defc}).  Then it represents the $n$-th component of the Chern character of the bundle
$E$.  
\end{theorem}
\begin{remark}
One can show that the second component of the Chern character -- the first Chern class
 -- can be described also by the following cocycle $\phi$, cohomologous
to $Ch^2(x_0,x_1,x_2)$.
Let the points $ x_0$, $x_1$ and $x_2$ be sufficiently close. $e(x_k)\in M_N(\mathbb{C})$
is an orthogonal projection and let $E_k\subset \mathbb{C}^N$
be its image, $k=0$, 1~, 2~.
Let $P_{kl}$ be the orthogonal projection from $E_k$ to $E_l$.
Then $P_{12}P_{23}P_{31}$ is a linear transformation in $E_1$, sufficiently close to 1,
and put
\[
\phi(x_0,x_1,x_2)=-\frac{1}{\pi}\Im \log \det(P_{12}P_{23}P_{31})
\]
where $\Im$ denotes the imaginary part.
 (Here $\log$ can be unambiguously defined by requiring that $\log(1)=0$, since
$\det(P_{12}P_{23}P_{31})$ is close to 1).
\end{remark}

To prove the theorem \ref{prf} we need several lemmas.  
\begin{lemma}\label{clos}
$Ch^n$ is an Alexander-Spanier cocycle.  
\end{lemma}
\begin{proof}  
First,  notice that $Ch^n$ is clearly a continuous function (in the neighborhood of the diagonal).  
Now suppose that we are given $n+2$ points $x_0,  \dots,  x_{n+1}$.
Similarly to (\ref{defe}) we can construct an idempotent over
$\Delta^{n+1}=\{t_j,  j=0,  \dots,  n+1 \mid \sum_{j=0}^{n+1} t_j =1,  t_j \geq 0\}$

\begin{equation}
F(t_0,  \dots,  t_{n+1})=\frac{1}{2\pi i}\int_{|\lambda-1|=1/2}
\left(\lambda-\sum_{j=0}^{n+1} t_je(x_j)\right)^{-1}d\lambda
\end{equation}

We can consider for any $n$ a form $TrF(dF)^n$ on $\Delta^{n+1}$.  This form is 0 for odd  $n$ and
is closed for even $n$.  Indeed,  since $F^2=F$,  we have $FdF +dFF=dF$.  Multiplying by $F$
we get $FdFF=0$.  Also $FdFdF=dFdF-dFFdF=dFdF-dFdF+dFdFF=dFdFF$.  Using this we get for $n$ odd:
\[
TrF(dF)^n=TrF(dF)^nF=TrFdFF(dF)^{n-1}=0
\]
For $n$ even we get:
\[
dTrF(dF)^n=dTrF(dF)^nF=Tr(dF)^{n+1}F+TrF(dF)^{n+1}=0,
\] 
where we have used the calculation for odd $n$.  
Now,  the restriction of $F$ to the face given by equation $t_j=0$ coincides (after an obvious
renumeration of variables) with the idempotents
constructed by the formula (\ref{defe}) from the points $x_0,  \dots,  \hat x_j,  \dots,  x_{n+1}$.  
This implies that
\begin{multline}
\partial Ch^n(x_0,  \dots,  x_{n+1})=\sum_{j=0}^{n+1} (-1)^jCh^n(x_0,  \dots,  \hat x_j,  \dots,  x_{n+1})\\
=b_n\int_{\partial \Delta^{n+1}}Tr F(dF)^n=
\int_{  \Delta^{n+1}} dTrF(dF)^n=0
\end{multline}
by  Stokes' theorem,$b_n=\frac{(-1)^{n/2}}{(2\pi i)^{n/2}(n/2)!}$.
\end{proof} 

\begin{lemma} \label{indp}
The cohomology class of $Ch^n$depends only on the bundle $E$,  and not on the particular
 idempotent  representing the bundle.  
\end{lemma}
\begin{proof}
Suppose that we have two different embeddings of our bundle
into trivial bundles,  and let $e_0(x)$,  $e_1(x)$ be the two corresponding idempotents.  We can suppose
that $e_0$ and $e_1$ are homotopic (this can be achieved by ``enlarging'' original trivial bundles).  
This homotopy can be chosen to be at least piecewise-smooth.  Now we will show that the cocycles
corresponding to the two smoothly homotopic idempotents differ by a coboundary.  Let $e_{\tau}(x)$
be this homotopy.  Define by (\ref{defe}) the idempotent over $\Delta^n$ $e_{\tau}(t_1,  \dots.  t_n)$.  In the proof
of this lemma we will write it just as $e_{\tau}$.
Notice that for each $\tau$ the restriction of $ e_{\tau}(t_1,  \dots,  t_n)$ to the face $t_j=0$
depends only on $e(x_k),  k\neq j$.  Let $Ch^n_\tau$ be defined similarly by (\ref{defc}).   

\begin{multline}
\frac{d}{d\tau}Ch^n= b_n \int \limits_{\sum_{j=1}^n t_j\leq 1,  t_j \geq 0}
\frac{d}{d\tau}Tr\ e_{\tau}(t_1,  t_2,  \dots,  t_n)(de_{\tau}(t_1,  t_2,  \dots,  t_n))^n
\end{multline}
$b_n=\frac{(-1)^{n/2}}{(2\pi i)^{n/2}(n/2)!}$.
But
\begin{multline}
\frac{d}{d\tau}Tr\ e_{\tau}(de_{\tau})^n
=Tr\ \frac{d}{d\tau}(e_{\tau})(de_{\tau})^n+
\sum_{j=1}^nTr\ e_{\tau}(de_{\tau})^{j-1}
(\frac{d}{d\tau}de_{\tau})(de_{\tau})^{n-j}\\
=d\left(\sum_{j=1}^n(-1)^jTr\ e_{\tau}(de_{\tau})^{j-1}
(\frac{d}{d\tau}e_{\tau})(de_{\tau})^{n-j}\right)-\sum_{j=1}^nTr\ (de_{\tau})^{j}
(\frac{d}{d\tau}e_{\tau})(de_{\tau})^{n-j}\\
=d\left(\sum_{j=1}^n(-1)^jTr\ e_{\tau}(de_{\tau})^{j-1}
(\frac{d}{d\tau}e_{\tau})(de_{\tau})^{n-j}\right)
\end{multline}
since $\sum_{j=1}^nTr\ (de_{\tau})^{j}
(\frac{d}{d\tau}e_{\tau})(de_{\tau})^{n-j}=0$.  

Indeed,  we have as before
\[
\frac{d}{d\tau} e_{\tau}=\frac{d}{d\tau} e_{\tau}^2=e_{\tau}\frac{d}{d\tau} e_{\tau}+
(\frac{d}{d\tau} e_{\tau}) e_{\tau}
\]
and from this $e_{\tau}(\frac{d}{d\tau} e_{\tau})e_{\tau}=0$.
Hence
\begin{multline}
\sum_{j=1}^nTr\ (de_{\tau})^{j}
(\frac{d}{d\tau}e_{\tau})(de_{\tau})^{n-j}\\
=\sum_{j=1}^nTr\ (de_{\tau})^{j}
(\frac{d}{d\tau}e_{\tau})e_{\tau}^2(de_{\tau})^{n-j}+\sum_{j=1}^nTr\ (de_{\tau})^{j}e_{\tau}^2
(\frac{d}{d\tau}e_{\tau})(de_{\tau})^{n-j}\\
=2\sum_{j=1}^nTr\ (de_{\tau})^{j}
\left(e_{\tau}(\frac{d}{d\tau}e_{\tau})e_{\tau}\right)(de_{\tau})^{n-j}=0
\end{multline}

 Notice that the form $ \sum_{j=1}^n(-1)^jTr\ e_{\tau}(de_{\tau})^{j-1}
(\frac{d}{d\tau}e_{\tau})(de_{\tau})^{n-j}$ restricted to the face of $\Delta^n$ depends only on the values
of $e_{\tau}$ at the vertices of this face.  This allows us to define 
an  Alexander-Spanier $n-1$ cochain $T_{\tau}(x_0,  \dots,  x_{n-1})$.  
Construct,  by (\ref{defe}),  $e_{\tau}=e_{\tau}(t_1,  \dots,  t_{n-1})$.  Then 
\begin{multline}
T_{\tau}(x_0,  \dots,  x_{n-1})\\
=\frac{(-1)^{n/2}}{(2\pi i)^{n/2}(n/2)!} \int \limits_{\sum_{j=1}^{n-1} t_j\leq 1,  t_j \geq 0}
 \sum_{j=1}^n(-1)^jTr\ e_{\tau}(de_{\tau})^{j-1}
(\frac{d}{d\tau}e_{\tau})(de_{\tau})^{n-j}
\end{multline}
Now, with $b_n=\frac{(-1)^{n/2}}{(2\pi i)^{n/2}(n/2)!}$
\begin{multline}
\frac{d}{d\tau}Ch^n_{\tau}= b_n \int \limits_{\sum_{j=1}^n t_j\leq 1,  t_j \geq 0}
\frac{d}{d\tau}Tr\ e_{\tau}(t_1,  t_2,  \dots,  t_n)(de_{\tau}(t_1,  t_2,  \dots,  t_n))^n\\
=b_n\int \limits_{\partial \Delta^n} \sum_{j=1}^n(-1)^jTr\ e_{\tau}(de_{\tau})^{j-1}
(\frac{d}{d\tau}e_{\tau})(de_{\tau})^{n-j}\\
=\sum_{j=0}^n(-1)^jT_{\tau}(x_0,  \dots,  \hat x_j,  \dots,  x_{n})\\
=(\partial T_{\tau})(x_0,  \dots,  x_n)
\end{multline}
by  Stokes' theorem,  and
\begin{equation}
Ch^n_1-Ch^n_0=\partial \int_{0}^{1}T_{\tau}d\tau
\end{equation}	
This proves the lemma.
\end{proof}
 
\begin{lemma}\label{smth}
If $E$ is a smooth vector bundle over a smooth manifold,  
 $Ch^n$ represents the $n$-th component of the Chern character.  
\end{lemma}
\begin{proof}
We will show that under the canonical map from the Alexander-Spanier
complex to de the Rham cohomological
complex $Ch^n$ is mapped
to the differential form representing the $n$-th component of Chern character constructed from the curvature
of the natural connection on $E$ induced by the embedding in the trivial bundle.  
Namely,  $e(x)$ defines connection by $\nabla \xi=ed\xi,  \xi \in \Gamma (E)$.  The curvature of this connection
is $\nabla^2=edede$.  By the Chern-Weil theory differential form $\frac{(-1)^{n/2}}{(2\pi i)^{n/2}(n/2)!}Tr\ e(de)^n$
represents the component of the Chern character of the bundle $E$ in the $n$-th cohomology.  
Now we will compute the image of $Ch^n$ under the isomorphism with de Rham cohomology.  
  We consider an arbitrary point $x_0$ and $n$ curves
$x_j(\epsilon_j),  j=1,  \dots,  n$ such that $x_j(0)=x_0$  and let the tangent vector
of $x_j$ at $x_0$ be $v_j$.
We compute ($\delta$ is from (\ref{defd})):
\[
 \delta(x_0,  x_1(\epsilon_1),  \dots,  x_n(\epsilon_n))|_{\epsilon_j= 0}=0
\]
\[
\frac{\partial}{\partial \epsilon_k}\delta|_{\epsilon_j=0}=t_k(\mathcal{L}_{v_k}e)(x_0)
\]
Here $\mathcal{L}_{v_k}$ denotes directional derivative.  Using this and differentiating (\ref{frme})
we get:
\begin{multline}
\frac{\partial}{\partial \epsilon_k}e(t_1,  \dots,  t_n)|_{\epsilon_j=0}\\=
(1-e(x_0))(\mathcal{L}_{v_k}e)(x_0)e(x_0)+e(x_0)(\mathcal{L}_{v_k}e)(x_0)(1-e(x_0))
=\mathcal{L}_{v_k}e(x_0)
\end{multline}
Here we used the identities $(\mathcal{L}_{v_k}e)(x_0)=e(x_0)(\mathcal{L}_{v_k}e)(x_0)+(\mathcal{L}_{v_k}e)(x_0)e(x_0)$
and $e(x_0)(\mathcal{L}_{v_k}e)(x_0)e(x_0)=0$ obtained,  as before,  by differentiating the relation $e(x)^2=e(x)$,  and
then multiplying it by $e(x)$.  From this and (\ref{frme}) we obtain:
\begin{align}
de(t_1,  \dots,  t_n)|_{\epsilon_j=0}=0\\
\frac{\partial}{\partial \epsilon_k}de(t_1,  \dots,  t_n)|_{\epsilon_j=0}=(\mathcal{L}_{v_k}e)(x_0)dt_k
\end{align}
We can now  compute ( with $b_n=\frac{(-1)^{n/2}}{(2\pi i)^{n/2}(n/2)!}$) :

\begin{multline}
\frac{\partial}{\partial \epsilon_1} \dots \frac{\partial}{\partial \epsilon_n}
Ch^n\left( x_0,  x_1(\epsilon_1),  \dots,  x_n(\epsilon_n) \right) |_{\epsilon_j=0}\\
=b_n\int\limits_{\sum_{j=1}^n t_j\leq 1,  t_j \geq 0}
Tr\ \frac{\partial}{\partial \epsilon_1} \dots \frac{\partial}{\partial \epsilon_n}
e(t_1,  t_2,  \dots,  t_n)(de(t_1,  t_2,  \dots,  t_n))^n|_{\epsilon_j=0}\\
=b_n\int\limits_{\sum_{j=1}^n t_j\leq 1,  t_j \geq 0}
\sum_{\sigma\in S_n} sgn(\sigma)e(x_o)(\mathcal{L}_{v_{\sigma 1}}e)(x_0)\dots(\mathcal{L}_{v_{\sigma n}}e)(x_0)
dt_1\dots dt_n\\
=b_n\frac{1}{n!}
\sum_{\sigma\in S_n} sgn(\sigma)e(x_o)(\mathcal{L}_{v_{\sigma 1}}e)(x_0)\dots(\mathcal{L}_{v_{\sigma n}}e)(x_0)\\
=b_n e(de)^n_{x_0}(v_1,  \dots,  v_n)
\end{multline}
and the result is already antisymmetric in $v_0, \dots, v_n$.
\end{proof}

We can now prove the main result of this section.  

\begin{proof}[Proof of Theorem \ref{prf}.   ]
 In the smooth situation this result is proved in Lemma \ref{smth}.  For the general case
 notice that $E$ is isomorphic to a pull-back of some smooth vector bundle over a smooth manifold,
given by the classifying map into some Grassmanian.
As our considerations are functorial under continuous mappings and, in the smooth
case, our cocycles represent the components of the Chern character the statement
of the theorem is true for some embedding of $E$ into a trivial bundle.   
  But then by Lemma \ref{indp} it is true for
 any such embedding.  
\end{proof}


\begin{thebibliography}{999}
\bibitem[CM]{[CM]}
A. Connes, H. Moscovici. Cyclic cohomology, the Novikov conjecture, and hyperbolic groups,
{\itshape Topology} 29:3 (1990), 345-388.

\bibitem[CST]{[CST]}
A. Connes, D. Sullivan, N.Teleman. Quasiconformal mappings, operators on Hilbert space
and local formulae for characteristic classes, {\itshape Topology }33 (1994),663-681.

\bibitem[Ge]{[Ge]}
E. Getzler. The odd Chern character in cyclic homology and spectral flow.
{\itshape Topology }25 (1993), 489-507

\bibitem[MW]{[MW]}
H. Moscovici, F.-B. Wu, Localization of topological Pontrjagin classes via finite
propagation speed, {\itshape Geom. and Func. Analysis} 4 (1994), 52-92.
\end{thebibliography}
\end{document}